\documentclass[12pt, reqno]{amsart}
\newtheorem{theorem}{Theorem}[section]
\newtheorem{lemma}[theorem]{Lemma}

\theoremstyle{definition}
\newtheorem{definition}[theorem]{Definition}

\theoremstyle{remark}
\newtheorem{remark}[theorem]{Remark}

\numberwithin{equation}{section}

\def\bpf{\begin{proof}}
\def\epf{\end{proof}}

\usepackage{hyperref}

\usepackage{amssymb,amsmath,graphicx,amsthm}

\def\simto{\overset\sim\to}

\def\simleq{\underset\sim<}

\def\simgeq{\underset\sim>}

\def\simle{\underset\sim<}

\def\simge{\underset\sim>}

\def\T{\text}

\def\1#1{\overline{#1}}

\def\2#1{\widetilde{#1}}

\def\3#1{\widehat{#1}}

\def\4#1{\mathbb{#1}}

\def\5#1{\frak{#1}}

\def\6#1{{\mathcal{#1}}}

\def\C{{\4C}}

\def\R{{\4R}}

\def\N{{\4N}}

\def\Z{{\4Z}}
\def\A{\6A}
\def\M{\6M}
\def\N{\6N}

\def\La{\Lambda}

\def\sumK{\underset{|K|=k-1}{{\sum}'}}

\def\sumJ{\underset{|J|=k}{{\sum}'}}

\def\T{\text}
\newcommand{\Om}{\Omega}
\newcommand{\om}{\omega}

\newcommand{\bom}{\bar{\omega}}

\newcommand{\we}{\wedge}

\newcommand{\no}[1]{\|{#1}\|}

\def\R{{\Bbb R}}

\def\C{{\Bbb C}}

\def\Z{{\Bbb Z}}

\def\la{\langle}
\def\ra{\rangle}
\def\di{\partial}
\def\dib{\bar\partial}
\def\Label#1{\label{#1}}

\def\simto{\overset\sim\to}

\def\simleq{\underset\sim<}

\def\simgeq{\underset\sim>}

\def\simle{\underset\sim<}

\def\simge{\underset\sim>}

\def\T{\text}

\def\1#1{\overline{#1}}

\def\2#1{\widetilde{#1}}

\def\3#1{\widehat{#1}}

\def\4#1{\mathbb{#1}}

\def\5#1{\frak{#1}}

\def\6#1{{\mathcal{#1}}}

\def\C{{\4C}}

\def\R{{\4R}}

\def\N{{\6N}}

\def\Z{{\4Z}}

\def\A{\6A}

\def\M{\6M}

\def\La{\Lambda}

\def\sumK{\underset{|K|=k-1}{{\sum}'}}

\def\sumJ{\underset{|J|=k}{{\sum}'}}

\numberwithin{equation}{section}

\def\T{\text}

\theoremstyle{plain}

\newtheorem{corollary}[theorem]{Corollary}

\newtheorem{proposition}[theorem]{Proposition}

\theoremstyle{definition}

\theoremstyle{remark}

\begin{document}

\title[Equivalence of estimates on a domain and its boundary]{Equivalence of estimates on a domain and its boundary}

\author{Tran Vu Khanh}
\address{T.V.~Khanh}

\address{Tan Tao University, Tan Duc e-city, Long An province, Vietnam}
\email{khanh.tran@ttu.edu.vn}

\address{School of Mathematics and Applied Statistics, University of Wollongong, NSW 2522 Australia}
\email{tkhanh@uow.edu.au}

\maketitle

\begin{abstract} Let $\Om$ be a pseudoconvex domain in $\C^n$ with smooth boundary $b\Om$. We define general estimates $(f\T-\M)^k_{\Om}$ and $(f\T-\M)^k_{b\Om}$ on $k$-forms for the complex Laplacian $\Box$ on $\Om$ and the Kohn-Laplacian $\Box_b$ on $b\Om$. For $1\le k\le n-2$, we show that $(f\T-\M)^k_{b\Om}$ holds if and only if $(f\T-\M)^k_{\Om}$ and $(f\T-\M)^{n-k-1}_{\Om}$ hold. Our proof relies on Kohn's method in \cite{Koh02}.  
\end{abstract}

%\section{Introduction}

\tableofcontents
\section{Introduction and results}
Let $\Om$ be a pseudoconvex domain $\C^n$ with smooth boundary $b\Om$. Let $L_2^{0,k}(\Om)$  be the space of square-integrable $(0,k)$-forms (or $k$-forms for short) on $\Om$. We have a  complex of densely defined operators $\bar\partial$ with $L^2$-adjoint $\bar\partial^*$ 
\begin{equation}\begin{split}
L_2^{0,k-1}(\Om)\underset{\dib^*}{\overset{\dib}{\rightleftarrows}} L_2^{0,k}(\Om)\underset{\dib^*}{\overset{\dib}{\rightleftarrows}} L_2^{0,k+1}(\Om),
 \end{split}\end{equation}
and the complex Laplacian is defined by $\Box:=\dib^*\dib+\dib\dib^*:L_2^{0,k}(\Om)\to L_2^{0,k}(\Om)$. The inverse is called the $\dib$-Neumann operator. We refer the reader to \cite{FK72, CS01, Str10, Zam08} for background on the $\dib$-Neumann problem. A general estimate for the complex Laplacian was introduced in \cite{Kha10} as follows.

For $z_o\in b\Om$, we choose local real coordinates $(a,r)\in \R^{2n-1}\times \R$ at $z_o$ and denote by $\xi$ the dual variables to the $a$'s. For a smooth, non-decreasing function $f$ with $\frac{f(t)}{t^{1/2}}$ decreasing, we denote by $f(\Lambda)$ the pseudodifferential operator of symbol $f\left((1+|\xi|^2)^{\frac{1}{2}}\right)$ which is defined by $f(\Lambda)u=\mathcal F^{-1}\left(f\left((1+|\xi|^2)^{\frac{1}{2}}\right)\mathcal F u\right)$ for $u\in C^\infty_c$, where $\mathcal F$ is the Fourier transform in $\R^{2n-1}$.  We also work with the function multiplier $\M$, that means,  $\M$ is a smooth function in $\bar\Om$. 

\begin{definition}\Label{fMk} Then  the $\bar{\partial}$-Neumann problem is said to satisfy the {\it  $(f\T-\M)_\Om^k$ estimate} at  $z_0\in b\Om$  if there exist positive constants $c$, $C_\M$   and a neighborhood $U$ of $z_0$ such that 
\begin{equation}\Label{gd}\begin{split}
(f\T-\M)_\Om^k\qquad \no{ f(\Lambda) \M u}^2 \le c(\no{\dib u}^2+\no{\dib^*u}^2+\no{u}^2)+C_{\mathcal M}\no{u}^2_{-1}, 
\end{split}\end{equation}
for all $u\in C^\infty_c(U\cap \bar{\Omega})^{k}\cap \T{Dom}(\bar{\partial}^*)$ and $k\ge 1$.
\end{definition}
%Here, when $\Om$ is $q$-pseudoconvex the matrix multiplier $\M$ acting a $k$-form $u=\sumJ u_J\bom_J$ is defined by
%\begin{eqnarray}\Label{Mu2}
%\begin{split}
%|\M u|^2 :=&\sumK\sum_{ij=1}^{n} \m_{ij} u_{iK}\bar u_{jK}}\\ 
%(\T{resp. }\quad
%\M u :=&\sumJ\sum_{j=1}^{n} \overline{\M}_j u_{J}\bom_j\we \bom_J\quad).\\
%\end{split}
%\end{eqnarray}
 
%This estimate is the most general to date in the $L^2$-estimates for the $\dib$-Neumann problem. 
In particular, for suitable  choice of $f$ and $\M$, the estimate $(f\T-\M)^k_\Om$ becomes a subelliptic estimate (see \cite{Koh64, Koh79, Cat83, Cat87,KZ11a}), a superlogarithmic estimate (see \cite{Koh02,KZ10,KZ12}), a compactness estimate and a weak compactness estimate (see \cite{Cat84a,FS98, FS01, Har11,KZ12c,McN02c,Str08}). \\

On a hypersurface $M$ in $\C^n$, the Cauchy-Riemann operator $\dib$ induces in a natural way the tangential Cauchy-Riemann operator $\dib_b$. The $\dib_b$ complex has played an important role in the study of boundary values of holomorphic functions and in the  problem of holomorphic extension \cite{KR65}. When $M$ is pseudoconvex, Shaw\cite{Sha85b} and Kohn\cite{Koh86} proved that $\dib_b$ has closed range. In recent time, this result has been extended to a  CR manifold $M$ of hypersurface type   in $\C^n$ \cite{Nic06, KN06, Bar13}.

Let $\dib^*_b$ be the $L^2$-adjoint of $\dib_b$ and $\Box_b=\dib_b\dib^*_b+\dib^*_b\dib_b$, the Kohn-Laplacian. Similarly as above, a general estimate for the Kohn-Laplacian operator holds at $z_o\in M$ for $k$-forms, $1\le k\le n-2$,  if there exist positive constants $c$, $C_\M$ and a neighborhood $U$ of $z_0$ such that 
\begin{equation}\Label{gh}\begin{split}
(f\T-\M)_M^k\qquad \no{ f(\Lambda)\M  u}_b^2 \le c(\no{\dib_b u}_b^2+\no{\dib_b^*u}_b^2+\no{u}^2_b)+C_{\mathcal M}\no{u}^2_{b,-1}, 
\end{split}\end{equation}
for all $u\in C^\infty_c(U\cap M)^{k}$. %On the top degrees,  $(f\T-\M)^0_M$ and 
%$(f\T-\M)^{n-1}_M$ are defined by 
%\begin{equation}\Label{gh0}(f\T-\M)^0_M  \qquad\no{ f(\Lambda)\M \zeta u}^2_b\le c\no{\dib_b u}_b^2+C_\M\no{u}_{b,-1}^2,\end{equation}
%for all $\zeta\in C_c^\infty(U\cap  M)$ and $u\in C^\infty(M)\cap \T{Ker}(\dib_b)^\perp$; and 
%\begin{equation}\Label{ghn-1}(f\T-\M)^{n-1}_M  \qquad\no{f(\Lambda)\M\zeta u}^2_b\le c\no{\dib_b^* u}_b^2+C_\M\no{u}_{b,-1}^2,\end{equation}
%for all $\zeta\in C_c^\infty(U\cap  M)$ and $u\in C^\infty(M)^{n-1}\cap \T{Ker}(\dib_b^*)^\perp$.

 When $f(t)=t^\epsilon$ and $\M=1$, \eqref{gh} is a subelliptic estimate. Subelliptic estimates for $\Box_b$ are well understood when $\Om$ is a strongly pseudoconvex domain \cite{Koh85}, or a pseudoconvex domain of finite type with comparable Levi  form \cite{Koe02}.  When $f(t)=1$ and $\M$ is an arbitrary constant, \eqref{gh} is called a compactnes estimate. Straube and Raich \cite{RS08} showed that a compactness estimate for $\Box_b$ holds when $M$ satisfies properties $(P_k)$ and $(P_{n-k-1})$.\\

The problem we address in this paper is  the equivalences of the estimates for the complex Laplacian $\Box$ on a pseudoconvex domain and the corresponding estimates for the Kohn-Laplacian $\Box_b$ on the boundary. In fact, we shall prove the following results. 

\begin{theorem}\Label{t1} Let $\Om$ be a pseudoconvex domain in $\C^n$ ($n\ge 2$) with smooth boundary $b\Om$, $\M$ be a smooth function on $\bar\Om$ and $f:\R^+\to\R^+$ be a smooth function such that $\frac{f(t)}{t^{1/2}}$ increasing. Then, for any $1\le k\le n-2$, the estimate $(f\T-\M)^k_{b\Om}$ for $\Box_b$ holds if and only if the estimates $(f\T-\M)^k_{\Om}$ and $(f\T-\M)^{n-k-1}_{\Om}$ for $\Box$ hold.\\

%Moreover, on the top degrees, the estimates $(f\T-\M)^0_{b\Om}$ and $(f\T-\M)^{n-1}_{b\Om}$ for $\Box_b$ are a consequence of $(f\T-\M)^1_{\Om}$ for $\Box$.
\end{theorem}

Our proof relies on Kohn's method in \cite{Koh02}. We also work with $\M$ is a vector function. 
In particular, let $M$ be a pseudoconvex hypersurface in $\C^n$. Denote $\Om^+$ and $\Om^-$ the pseudoconvex and pseudoconcave side of $M$, respectively. Let $u=\sumK u_J\bom_J$ be a $(0,k)$-form and $\M=(\M_j)_{j=1}^n$ be a smooth vector on $\bar\Om$. we use notation  $\M u$ on pseudoconvex domain $\Om^+$ by 
$$\M u=\sumK\sum_{j=1}^n \M_{j}u_{jK}\bom_K,$$ 
and on pseudoconcave domain $\Om^-$ by
$$\M u=\sumJ\sum_{j=1}^n\overline\M_ju_J\bom_j\we \bom_J.$$ 
We use notation $(f\T-\M)^k_{\Om^\pm}$ in the obvious sense. On the hypersurface $M$, let $(f\T-\M)^k_{M,+}$ and $(f\T-\M)^k_{b\Om,-}$ denote general estimates acting on $(0,k)$-forms of positive Kohn's microlocalization and negative Kohn's microlocalization (see Section~2 below), respectively. We have the following equivalences
\begin{theorem}\Label{t2}Let $M$ be a pseudoconvex hypersurface in $\C^n$ and let $\M$ be a function/vector multiplier. Then we have
$$(f\T-\M)^k_{\Om^+}\Longleftrightarrow(f\T-\M)^k_{b\Om,+} \Longleftrightarrow (f\T-\M)^{n-k-1}_{b\Om,-}\Longleftrightarrow(f\T-\M)^{n-k-1}_{\Om^-},$$
for any $1\le k\le n-1$. 
\end{theorem}
The proof of Theorem~\ref{t2} is a consequence of the results in Section 5 and 6. The proof of Theorem~\ref{t1} follows immediately from Theorem~\ref{t2} since the elliptic estimate holds for $u^0$.

The paper will be presented as follows. In Section 2, we give a brief introduction to the Cauchy-Riemann operator, the tangential Cauchy-Riemann operators and Kohn microlocalization. The microlocal estimates on $b\Om$ and $\Om^\pm$ are given in Section 3 and 4. The main part of the proof for Theorem~\ref{t2} lies in Section 5. In Section 6, we give equivalence of microlocal estimates on hypersurface and complete the proof of Theorem \ref{t1}.

\section{Preliminaries}

Let $M$ be a smooth hypersurface in $\C^n$. We start by denoting  by  $\A_b^{0,k}$ the space of smooth sections of the vector bundle $\big (T^{0,1}(M)^*\big)^k$ on $M$.  The tangential Cauchy-Riemann operator $\dib_b : \A_b^{0,k}\to \A_b^{0,k+1} $ is defined as follows. If $u \in \A_b^{0,k}$,  let $u'$ be a (0,k)-form whose restriction to $M$  equals $u$. Then $\dib_b u$ is the restriction of the Cauchy-Riemann operator $\dib$ on $u'$ to $M$. We can define a Hermitian inner product on $\A^{0,k}_b$ by 
$$(u,v)_b=\int_M\la u,v \ra dS,$$    
where $dS$ is the surface element on $M$. The inner product gives rise to an $L_2$-norm $\no{\cdot}_b$.
 We define $\dib^*_b$ to be the $L^2$-adjoint of $\dib_b$ in the standard way. Thus $\dib^*_b:\A_b^{0,k+1}\to \A_b^{0,k} $  for $k\ge 0$. 
 For $u,\,v\in C^\infty_c(U\cap M)^k$, we denote the tangential energy by
$$Q_b(u,v)=(\dib_b u,\dib_b u)_b+(\dib^*_b u,\dib^*_bv)_b+(u,v)_b.$$

Let $z_0\in M$ and $U$ be a neighborhood of $z_0$; we fix a defining function $r$ of $M$ such that $|\di r|=1$ on $U\cap M$. We choose an orthonormal basis of $(1,0)$ forms $\om_1,\dots,\om_n=\di r$ and the dual basis of $(1,0)$ vector fields $L_1,...,L_n$ such that $L_j|_{z_0}=\di_{z_j}$; thus $L_1,\dots,L_{n-1}\in T^{1,0}M$ and $L_n(r)=1$. We define $T=\frac{1}{2}(L_n-\bar L_n)$ and $D_r=\frac{1}{2}(L_n+\bar L_n)$. It follows 
\begin{eqnarray}\Label{4.1a}
\begin{split}
L_n=D_r+T \T{~~~~and~~~~ }\bar L_n=D_r-T.
\end{split}
\end{eqnarray}

Denote $C_c^\infty(U\cap M)^k$ the forms of $\A^{0,k}_b$ with  compact support in $U$. 
We write $k$-forms in $C^\infty_c(U\cap M)^k$ as \begin{eqnarray}\Label{ubb}
u=\sumJ u_J\bom_J,
\end{eqnarray}
 where $J=\{j_1,\dots,j_k\}$ is a ordered multiindex and $\sum'$ denotes summation over strictly increasing index sets. If $J$ decomposes as $J=iK$, then $u_{iK}=sign\binom{J}{iK}u_J$. 
Then on $M$, the operator $\dib_b$ and $\dib^*_b$ are expressed by
\begin{equation}\Label{dib}
\dib_b u=\sumJ\sum_{j=1}^{n-1} \bar L_j u_J\bom_j\we\bom_J+...
 \end{equation}
and 
\begin{eqnarray}\Label{dib*}
\dib_b^*u=-{\sumK}\sum_j^{n-1} L_ju_{jK}\bom_K+... 
\end{eqnarray}
where dots refer the error terms in which $u$ is not differentiated. By developing the equalities \eqref{dib} and \eqref{dib*}, the key technical result is contained in the following.

\begin{proposition} \Label{l4.1}For two indices  $q_1,q_2$ ; $(0\le q_1\le q_2\le n-1)$, we have 
\begin{eqnarray}\begin{split}\Label{4.3}
Q_b(u,u):=&\no{\dib_b u}_b^2+\no{\dib^*_b u}_b^2+ \no{ u}_b^2\\
\simgeq& \sumK\sum_{ij=1}^{n-1}(r_{ij}Tu_{iK},u_{jK})_b-\sumJ\sum_{j=q_1}^{q_2}(r_{jj}Tu_J,u_J)_b\\
&+\frac{1}{2}\sumJ\Big(\sum_{j=1}^{q_1-1}\no{\bar L_j u_J}^2_b+\sum_{j=q_2+1}^{n-1}\no{\bar L_j u_J}^2_b+\sum_{j=q_1}^{q_2}\no{ L_j u_J}^2_b\Big).
\end{split}
\end{eqnarray}
for any $u\in C^\infty_c(U\cap M)^k$.
\end{proposition}
Here and in what follows, $\simge$ or $\lesssim$ denote inequality up to a constant; and $r_{ij}$ is the coefficient of $\om_i\we\bom_j$ when $\di\dib r$ expressed in this basis. We refer \cite{Kha10b} for the proof of this proposition. Note that, conversely, we have
\begin{eqnarray}\begin{split}\Label{4.4}
\no{\dib_b u}_b^2+&\no{\dib^*_b u}_b^2\\
\simleq& \left|\sumK\sum_{ij=1}^{n-1}(r_{ij}Tu_{iK},u_{jK})_b-\sumJ\sum_{j=q_1}^{q_2}(r_{jj}Tu_J,u_J)_b\right|\\
&+\sumJ\Big(\sum_{j=1}^{q_1-1}\no{\bar L_j u_J}^2_b+\sum_{j=q_2+1}^{n-1}\no{\bar L_j u_J}^2_b+\sum_{j=q_1}^{q_2}\no{ L_j u_J}^2_b\Big)+\no{u}_b^2,
\end{split}
\end{eqnarray}
 for all $u\in C_c^\infty(U\cap M)^k$ for any $k$.\\

In $U$, we choose special boundary coordinate $(x_1,...,x_{2n-1},r)$. Let $\xi=$ $(\xi_1,...,\xi_{2n-1})$ $=(\xi',\xi_{2n-1})$ be the dual coordinates to $\{x_1,...,x_{2n-1}\}$. 
We also decompose $(x_1,...,x_{2n-1})$ $=(x',x_{2n-1})$ so that $T^\C_{z_0}M$ is defined by $x_{2n-1}=0$ in $T_{z_0}M$.
Let $\psi^++ \psi^-+ \psi^0=1$ be a $C^\infty$ partition of the unity in the sphere $|\xi|=1$ such that
$\psi^\pm$ are 1 at the poles $(0,...,\pm 1)$ and $\psi^0=1$ at the equator, that is, at  $\xi_{2n-1}=0$. We extend these functions to $\R^{2n-1}\setminus\{0\}$ as homogeneous functions of degree 0. We may assume that the supports of the functions $\psi^+$, $\psi^-$ and $\psi^0$ are contained in the cones 
 \begin{eqnarray}\begin{split}
\mathcal C^+=&\{\xi \big |\xi_{2n-1}> \frac{1}{2} |\xi'| \};\\
\mathcal C^-=&\{\xi\big|-\xi_{2n-1}> \frac{1}{2} |\xi'|  \};\\
\mathcal C^0=&\{\xi \big | | \xi_{2n-1}|<|\xi'| \}.
\end{split}
\end{eqnarray}
Then $\T{supp}{\psi^+}\subset\subset  \mathcal C^+$, $\T{supp}{\psi^-}\subset\subset  \mathcal C^-$ and  $\T{supp}{\psi^0}\subset\subset \mathcal C^0$.\\

The operators $\Psi=\Psi^{\overset\pm0}$ with symbols $\psi=\psi^{\overset\pm0}$ are defined by 
$$\widetilde{\Psi \varphi}(\xi)=\psi(\xi)\tilde{\varphi}(\xi) \qquad \T{for }\quad \varphi\in C^\infty_c(U\cap M);$$
$$\widetilde{\Psi \varphi}(\xi,r)=\psi(\xi)\tilde{\varphi}(\xi,r) \qquad \T{for }\quad \varphi\in C^\infty_c(U\cap \Om).$$
The microlocal decomposition  $\varphi=\varphi^++\varphi^-+\varphi^0$ of a function $\varphi\in C^\infty_c(U\cap M)$ is defined by 
$$\varphi=\zeta \Psi^+\varphi+\zeta \Psi^-\varphi+\zeta \Psi^0\varphi,$$
 where $\zeta\in C^\infty(U'), \bar U\subset \subset U'$ and $\zeta\equiv1$ on $U$.  \\
 For a form $u$, the microlocal decomposition $u=u^++u^-+u^0$ is accordingly defined coefficientwise. 
We will define $(f\T-\M)^k_{M,+}$ and $(f\T-\M)^k_{M,-}$ which were mentioned in Section 1.\\

Let $u$ be as in \eqref{ubb}. A {\it function multiplier} $\M\in \A_b^{0,0}$ on the positive (resp. negative) microlocalization $u^+$ (resp. $u^-$)  is defined by $$\M u^+=\sumJ\M u^+_J\bom_J \quad (\T{resp.} \quad \M u^-=\sumJ\overline{\M} u^-_J\bom_J),$$
and for full $u$ defined by $\M u=\sumJ\M u_J\bom_J$.\\

A $(1,0)$-form multiplier $\M\in\A^{1,0}_b$ on the positive (resp. negative) component of $u$ on $M$ is defined by 
\begin{eqnarray}\label{Mpm}
\begin{split}
\M u^+ :=&\sumK\sum_{j=1}^{n-1} \M_j u^+_{jK}\bom_K\\ 
(\T{resp. }\quad
\M u^- :=&\sumJ\sum_{j=1}^{n-1} \overline{\M}_j u^-_{J}\bom_j\we \bom_J\quad).
\end{split}
\end{eqnarray}
\begin{definition} Let $M$ be a hypersurface, $z_0\in M$, and $\M\in \A_b^{0,0}$ or $\A^{1,0}_b$. we say that a $(f\T-\M)^k_{M,+}$- (resp. $(f\T-\M)^k_{M,-}$-) estimate  for $(\dib_b,\dib_b^*)$  holds at $z_0$ if $(f\T-\M)^k_{M}$ holds with $u$ replaced by $u^+$ (resp.  $u^-$), that is,  
$$(f\T-\M)^k_{M,+}\qquad \no{f(\Lambda)\M u^+}^2_b\lesssim Q_b(u^+,u^+)+\no{\Psi^+u}^2_{b,-1}+C_\M\no{u^+}_{b,-1}^2$$ 
(resp.
$$(f\T-\M)^k_{M,-}\qquad \no{f(\Lambda)\M u^-}^2_b\lesssim Q_b(u^-,u^-)+\no{\Psi^-u}^2_{b,-1}+C_\M\no{u^-}_{b,-1}^2\quad ).$$ 
 \end{definition}

The hypersurface $M$ is said to be pseudoconvex at $z_0$ if either of the two components of $\C^n\setminus M$  is pseudoconvex at $z_0$. Denote by $\Om^+=\{z\in U | r(z)<0 \}$ the pseudoconvex side of $M$ and by $\Om^-$ the other one. Then $\Om^-$ is pseudoconcave. We also use the notation $\omega_n^\pm=\pm\partial r$ for the exterior conormal vectors to $\Om^\pm$.

Under choice of such basis, we check readily that $u\in \T{Dom}(\dib^*)$ if and only if $u_{nK}|_{b\Om=M}=0$ for any $K$. For $u,\,v\in C^\infty_c(U\cap \bar\Om^\pm)^k\cap \T{Dom}(\dib^*)$, we denote the energy by
$$Q(u,v)=(\dib u,\dib v)+(\dib^* u,\dib^* v)+(u,v).$$
Using integration by parts for $u\in C^\infty_c(U\cap \bar\Om^\pm)^k\cap \T{Dom}(\dib^*)$, we have
\begin{eqnarray}\Label{kmh}\begin{split}
Q(u,u)\simge &\sumK\sum_{ij=1}^{n-1}\int_{b\Om}r_{ij}u_{iK}\bar u_{jK}dS-\sumJ\sum^{q_o}_{j=1}\int_{b\Om}r_{jj}|u_{J}|^2dS\\
&+\frac{1}{2}\left(\sum^{q_o}_{j=1}\no{L_j u}^2+\sum_{q_o+1}^{n}\no{\bar L_j u}^2 \right).
\end{split}
\end{eqnarray}
for any $0\le q_o\le n-1$. Also, we have the converse inequality 
\begin{eqnarray}\Label{new1}\begin{split}
Q(u,u)\lesssim  &\left|\sumK\sum_{ij=1}^{n-1}\int_{b\Om}r_{ij}u_{iK}\bar u_{jK}dS-\sumJ\sum^{q_o}_{j=1}\int_{b\Om}r_{jj}|u_{J}|^2dS\right|\\
&+\left(\sum^{q_o}_{j=1}\no{L_j u}^2+\sum_{q_o+1}^{n}\no{\bar L_j u}^2 \right)+\no{u}^2.
\end{split}
\end{eqnarray}
We finish this section with an estimate in the normal vector field $D_r$. 

\begin{lemma}\Label{l3.38}
If the $(f\T-\M)^k$ estimate holds, then we have
$$\no{\Lambda^{-1} D_r u}^2+\no{\Lambda^{-1}f(\Lambda)D_r\M u}^2\lesssim Q(u,u)+C_\M \no{u}_{-1}^2,$$
for any $u\in C^\infty_c(U\cap \bar\Om)^k\cap \T{Dom}(\dib^*)$. 
\end{lemma}
{\it Proof. }
Using \eqref{4.1a}, we have 
\begin{eqnarray}
\begin{split}
\no{\Lambda^{-1} D_r u}^2+\no{\Lambda^{-1}f(\Lambda)D_r\M u}^2\lesssim &\no{\La^{-1}\bar L_nu}^2+\no{\La^{-1}Tu}^2\\
&+\no{\Lambda^{-1}f(\Lambda)\bar L_n (\M u)}^2+\no{\Lambda^{-1}f(\Lambda)T \M u}^2\\
\lesssim &\no{\bar L_n u}^2+\no{f(\Lambda) \M u}^2+\no{u}^2\\
\lesssim& Q(u,u)+C_\M\no{u}_{-1}^2.
\end{split}
\end{eqnarray}
$\hfill\Box$\\

\section{Basic microlocal  estimates on $b\Om$}
In this section we prove the basic microlocal estimates on hypersurface. We first start with 0-Kohn microlocalization $u^0$.

\begin{lemma}  \Label{l4.4}
 Let $M$ be a  hypersurface and $z_0$ a point of $M$. Then there is a  neighborhood $U$ of $z_0$ such that
  $$ Q_b(u^0,u^0) \cong \no{ u^0}_{b,1}^2, $$
for all $u\in C_c^\infty(U\cap M)^k$ with any $k$. 
\end{lemma}
{\it Proof.} Using twice the inequality \eqref{4.3}  for $q_1=q_2=0$ and $q_1=0, q_2=n-1$ and taking summation, we get
\begin{eqnarray}\Label{4.5}
\begin{split}
Q_b(u,u)
\simgeq& \sumK\sum_{ij=1}^{n-1}(r_{ij}Tu_{iK},u_{jK})_b-\sumJ\sum_{j=1}^{n-1}(r_{jj}Tu_J,u_J)_b\\
&+\frac{1}{2}\Big(\sumJ\sum_{j=1}^{n-1}\no{L_j u_J}^2_b+\sumJ\sum_{j=1}^{n-1}\no{\bar L_j u_J}^2_b\Big)\\
\simgeq& \no{\Lambda' u}^2_b-\big(\epsilon+\T{diam}(U)\big)\no{Tu}^2_b-C_\epsilon \no{u}^2_b,
\end{split}
\end{eqnarray}
where $\La'$ is the pseudodifferential operator of order 1 whose symbol is  $(1+\sum_{j=1}^{2n-2}|\xi_j|^2)^\frac{1}{2}$. 
To explain the last estimate in \eqref{4.5}, recall that $L_j|_{z_0}=\partial_{z_j},\,\,j=1,...,n-1$ and that the coefficients of the $L_j$'s are $C^1$;  therefore, the third line of \eqref{4.5} is bounded below by $\no{\Lambda'u}^2-\T{diam}(U)\no{Tu}^2_b-\no{u}_b^2$.
Apply \eqref{4.5} for $u^0$ and notice that $\no{\Lambda'u^0}^2_b\simgeq \no{Tu^0}^2_b$. Taking $U$ and $\epsilon$ suitably small, we conclude
$$Q_b(u^0,u^0)\simge \no{\La' u^0}^2_b\simge \no{\La u^0}^2_b.$$ 
On the other hand, the converse inequality is always true. 

$\hfill\Box$

We now give the basis estimates for $u^+$ and $u^-$.

\begin{lemma}  \Label{l4.5}
Let $M$ be a pseudoconvex hypersurface at $z_0$. Then, for a neighborhood $U$ of $z_0$, and for $\zeta'\equiv1$ over $\T{supp}(\zeta)$ and $\psi^{\pm\,\prime}\equiv1$ on $\T{supp}\,\psi^\pm$, we have 
\begin{eqnarray}{\it (i) } \begin{split}  
 Q_b(u^+,u^+) &+\no{\Psi^+u}_{b,-1}^2\\
\cong& \sumK\sum_{ij=1}^{n-1}(r_{ij}\zeta' R^+ u^+_{iK},\zeta'R^+u^+_{jK})_b\\
&+\sumJ\sum_{j=1}^{n-1}\no{\bar L_j u^+_J}^2_b+\no{u^+}_b^2+\no{\Psi^+u}_{b,-1}^2,
\end{split}
\end{eqnarray}
for all $u\in C_c^\infty(U\cap M)^k$ with  $k\ge 1$, where $R^+$ is the pseudodifferential operator of order $\frac{1}{2}$ whose symbol is $\xi_{2n-1}^{\frac{1}{2}}\psi^{+\,\prime}(\xi)$. Similarly, we have 
 \begin{eqnarray}{\it  (ii) }\begin{split}  
 Q_b(u^-,u^-)&+\no{\Psi^-u}_{b,-1}^2\\
 \cong&\sumJ\sum_{j=1}^{n-1}(r_{jj}\zeta' R^-u^-_J,\zeta' R^-u^-_J)_b- \sumK\sum_{ij=1}^{n-1}(r_{ij} \zeta' R^- u^-_{iK},\zeta' R^-u^-_{jK})_b\\
&+\sumJ\sum_{j=1}^{n-1}\no{ L_j u^-_J}^2_b+\no{u^-}_b^2+\no{\Psi^-u}_{b,-1}^2,
\end{split}
\end{eqnarray}
for any  $u\in C_c^\infty(U\cap M)^k$ with  $k\le n-2$, where $ R^-$ is the pseudodifferential operator of order $\frac{1}{2}$ whose symbol is $(-\xi_{2n-1})^{\frac{1}{2}}\psi^{-\,\prime}(\xi)$.
\end{lemma}
\begin{remark}\Label{convex-concave}
 Since $M$ is pseudoconvex at $z_0$, then there is a defining function of $M$ which satisfies on $M$
\begin{eqnarray}\Label{convex}\sumK\sum_{ij=1}^{n-1}r_{ij}u_{iK}\bar u_{jK}\ge 0,\end{eqnarray}
for any $u\in C_c^\infty(U\cap M)^k$ with $k\ge 1$, and also satisfies on $M$
\begin{eqnarray}\Label{concave}\sum_{j=1}^{n-1}r_{jj}|u|^2-\sumK\sum_{ij=1}^{n-1}r_{ij}u_{iK}\bar u_{jK}\ge 0,\end{eqnarray}
for any $u\in C_c^\infty(U\cap M)^k$ with $k\le n-2$,
 where $U$ is a neighborhood of $z_0$.
\end{remark}
{\it Proof.}  {\it (i). } 
We have
$$u_{iK}^+=\zeta \Psi^{+}u_{iK}=\zeta (\Psi^{+\,\prime})^2 \Psi^+u_{iK} = (\Psi^{+\,\prime})^2\zeta \Psi^+u_{iK}+[\zeta, (\Psi^{+\,\prime})^2]\Psi^+u_{iK}.$$
Since the supports of symbols of $\Psi^+$ and $[\zeta, (\Psi^{+\,\prime})^2]$ are disjoint, the operator $[\zeta, (   \Psi^{+\,\prime})^2]\Psi^+$ is of order $-\infty$ and we have
\begin{eqnarray}\Label{4.8}\begin{split}
(r_{ij}Tu_{iK}^+, u_{jK}^+)_b=&(r_{ij}T\zeta \Psi^+u_{iK}, \zeta \Psi^+u_{jK})_b\\
&=(r_{ij}T (\Psi^{+\,\prime})^2\zeta \Psi^+u_{iK}, \zeta \Psi^+u_{jK})_b+O(\no{\Psi^+u}^2_{b,-1})\\
&=((\zeta')^2 r_{ij}R^{+\,*}R^+\zeta \Psi^+u_{iK}, \zeta \Psi^+u_{jK})_b+O(\no{\Psi^+u}^2_{b,-1})\\
&=(r_{ij} \zeta' R^+\zeta \Psi^+u_{iK},\zeta' R^+\zeta \Psi^+u_{jK})_b\\
&+([(\zeta')^2 r_{ij},R^{+\,*}]R^+\zeta \Psi^+u_{iK}, \zeta \Psi^+u_{jK})_b+O(\no{\Psi^+u}^2_{b,-1}).
\end{split}\end{eqnarray}
From the pseudodifferential operator calculus we get
 \begin{eqnarray}\begin{split}
([\zeta^{\prime2} r_{ij},R^{+*}]R^+\zeta \Psi^+u_{iK}, \zeta \Psi^+u_{jK})_b|\lesssim \no{u^+}^2_b.
\end{split}\end{eqnarray}
Thus,
\begin{eqnarray}\Label{4.10}\begin{split}
&\sumK\sum_{ij=1}^{n-1}(r_{ij}T u^+_{iK},u^+_{jK})_b\\
=& \sumK\sum_{ij=1}^{n-1}(r_{ij}\zeta' R^+ u^+_{iK},\zeta'R^+u^+_{jK})_b+O(\no{u^+}^2_b)+O(\no{\Psi^+u}^2_{b,-1}).
\end{split}\end{eqnarray}
By Remark \ref{convex-concave}, the first sum in second line in \eqref{4.10} is nonnegative if $k\ge 1$.  Thus the first part of Lemma \ref{l4.5}  is proven  by applying Proposition \ref{l4.1} to $u^+$ with $q_1=q_2=0$.\\

{\it (ii). } The proof of the second part is similar. We have to notice that
\begin{eqnarray}\begin{split}
&\sumK\sum_{ij=1}^{n-1}(r_{ij}T u^-_{iK},u^-_{jK})_b-\sumJ\sum_{j=1}^{n-1}(r_{jj}Tu^-_J, u^-_J)_b\\
=&-\sumK\sum_{ij=1}^{n-1}(r_{ij}\zeta'  R^- u^-_{iK},\zeta' R^-u^-_{jK})_b+\sumJ\sum_{j=1}^{n-1}(r_{jj}\zeta' R^-u^-_J,\zeta' R^-u^-_J)_b\\
&+O(\no{u^-}^2_b)+O(\no{\Psi^-u}^2_{b,-1}),
\end{split}\end{eqnarray}
By Remark \eqref{convex-concave}, the second line is nonnegative for any $k$-form $u$ with $k\le n-2$. Thus the second part of Lemma \ref{l4.5}  is proven  by applying Proposition \ref{l4.1} to $u^-$ with $q_1=1$ and $q_2=n-1$.\\

$\hfill\Box$\\

\section{Basic  microlocal estimates on $\Om^+$ and $\Om^-$}

In this section,  we prove  the basic microlocal estimates on $\Om^+$ and $\Om^-$. We begin by introducing  the harmonic extension of a form from $b\Om$ to $\Om$ following  Kohn  \cite{Koh86} and \cite{Koh02}. \\

In terms of special boundary coordinate $(x,r)$, the operator $L_j$ can be written as
$$L_j=\delta_{jn}\frac{\di}{\di r}+\sum_k a_j^k(x,r) \frac{\di}{\di x_k}$$
for $j=1,. . ., n$. We define the tangential symbols of $L_j$, $1\le j\le n-1$, by
$$\sigma(L_j)((x,r),\xi)=-i \sum_k a_j^k(x,r) \xi_k,$$ 
 and 
$$\sigma(T)((x,r),\xi)=\frac{-i}{2} \sum_k \big(a_n^k(x,r)-\bar a_n^k(x,r)\big) \xi_k.$$
Note that $\sigma(T)$ is real. We set $$\sigma_b(L_j)(x,\xi)=\sigma({L_j})((x,0),\xi) \T{~~and~~ } \sigma_b(T)(x,\xi)=\sigma(T)((x,0),\xi)$$ 
and $$\mu(x,\xi)=\sqrt{\sum_j|\sigma_b(L_j)(x,\xi)|^2+ |\sigma_b(T)(x,\xi)|^2+1}.$$
Remember the notation $\Lambda_\xi=(1+|\xi|^2)^{\frac12}$;
in a neighborhood of $z_0$, we have $\mu(x,\xi)\sim\Lambda^1_\xi$. 
\\

Harmonic extension of boundary functions is defined as follows. Let $\varphi\in C^\infty(b\Om)$; define $\varphi^{(h)}\in C^\infty(\bar\Om^{\pm})$ by 
$$\varphi^{(h)}(z)=\sum_\nu \varphi^{(h)}_\nu(x_\nu,r)$$
where 
  $$\varphi_\nu^{(h)} (x_\nu, r)=(2\pi)^{-2n+1}\int_{\R^{2n-1}}e^{i x_\nu\cdot\xi_\nu} e^{\pm r\mu(x_\nu,\xi_\nu)}\widetilde{(\zeta_\nu\varphi^\nu)}(\xi_\nu)d\xi_\nu,$$
so that $\varphi^{(h)}=\varphi$ on $b\Om$ and therefore $\Psi^{\pm}(\varphi-\varphi^{(h)})=0$ on $b\Om$. Here $\{\zeta_\nu\}_\nu$ is a partition of unity subordinate to the covering $\{U_\nu\}_\nu$ of the boundary satisfying $\sum_\nu \zeta_\nu=1$, and $\phi^\nu$ is the function expressed in the local coordinates $(x_\nu,0)$ on $U_\nu$.

This extension is called  ``harmonic"  since $\triangle \varphi^{(h)} (x, r)$ has order 1 on $M$. In fact, we have 
\begin{eqnarray}\begin{split}
\triangle=&-\sum_{j=1}^n \frac{\di^2}{\di z_j \di \bar z_j}\\
=&-\sum_{j=1}^n L_j\bar L_j+ \sum_{k=1}^{2n-1}a^k(x,r)\frac{\di}{\di x_k}+a(x,r)\frac{\di}{\di r}\\
=&-\frac{\di^2}{\di^2 r}+T^2-\sum_{j=1}^{n-1} L_j\bar L_j+ \sum_{k=1}^{2n-1}b^k(x,r)\frac{\di}{\di x_k}+b(x,r)\frac{\di}{\di r}
\end{split}
\end{eqnarray}
since \eqref{4.1a} implies that $L_n\bar L_n=\frac{\di^2}{\di^2 r}-T^2+D$, where $D$ is a first order operator. Hence if $(x,r)\in U\cap \bar\Om^+, $
 \begin{eqnarray}\begin{split}\Label{4.15b}
\triangle (\varphi^{(h)})(x,r)=(2\pi)^{-2n+1}\int e^{ix\cdot \xi}e^{r\mu(x,\xi)}\big(p^1(x,r,\xi)+rp^2(x,r,\xi)\tilde\varphi(\xi,0)\big)d\xi,
\end{split}
\end{eqnarray}
where $p^k(x,r,\xi)$ denotes a symbol of order $k$, uniformly in $r$. For future use, we prepare the notation $P^1+rP^2$ for the pseudodifferential operator with symbol $p^1+rp^2$ which appears in the right of \eqref{4.15b}.
Along with \eqref{4.15b} we have 
$$L_j\varphi^{(h)}(x,r)=(L_j\varphi)^h(x,r)+E_j\varphi(x,r)$$
where 
$$E_j\varphi (x,r)= (2\pi)^{-2n+1}\int e^{ix\cdot \xi}e^{r\mu(x,\xi)}\big(p_j^0(x,r,\xi)+rp_j^1(x,r,\xi)\tilde\varphi(\xi)\big)d\xi$$
and 
$$\bar L_j\varphi^{(h)}(x,r)=(\bar L_j\varphi)^h(x,r)+\bar E_j\varphi(x,r)$$
for $j=1,...,n-1$.

\begin{lemma}\Label{l4.10}
For any $k\in \Z$ with $k\ge 0$, $s\in \R$, we have
\begin{enumerate}
  \item[(i)] $|||r^kf(\Lambda)\varphi^{(h)}|||_s\simle \no{f(\Lambda)\varphi}_{b, s-k-\frac{1}{2}}$,
  \item[(ii)] $|||D_rf(\Lambda)\varphi^{(h) }|||_s\simle \no{f(\Lambda)\varphi}_{b,s+\frac{1}{2}}$
  \end{enumerate}
for any $\varphi\in C^\infty_c(U\cap \bar\Om^+).$
\end{lemma}
{\it Proof. } We notice again that $\mu(x,\xi)\cong (1+|\xi|^2)^{\frac{1}{2}}$ over a small neighborhood of $z_0$, and then the proof of this lemma is similar to the proof  of Lemma 8.4 in \cite{Koh02}.\\

$\hfill\Box$\\

We define $\varphi_b$ to be the restriction of $\varphi\in C^\infty_c( U\cap \bar\Om^+)$ to the boundary. We have the elementary estimate
\begin{eqnarray}\Label{4.16b}
\no{\varphi_b}_{b, s}^2\lesssim |||\varphi|||_{s+\frac{1}{2}}^2+|||D_r\varphi |||^2_{s-\frac{1}{2}}.
\end{eqnarray}

The following lemma states the basic microlocal estimates on $\Om^+$.
\begin{lemma}\Label{l4.6} Let $\Om^+$ be pseudoconvex  at  $z_0$. If U is  a sufficiently small  neighborhood of $z_0$, then
we have the three estimates which follow
 \begin{eqnarray}\Label{4.12}|||\Psi^0\varphi|||^2_{1}\lesssim \sum_{j=1}^{n}\no{\bar L_j\Psi^0\varphi}^2+\no{\Psi^0\varphi}^2\quad
 \T{for any } \varphi \in C_c^\infty(U\cap \bar\Om^+),
 \end{eqnarray}
\begin{eqnarray}\Label{4.13}|||\Psi^-\varphi|||^2_{1}\lesssim \sum_{j=1}^{n}\no{\bar L_j\Psi^-\varphi}^2+\no{\Psi^-\varphi}^2\quad
\T{for any } \varphi \in C_c^\infty(U\cap \bar\Om^+),
\end{eqnarray}\begin{eqnarray}\Label{4.14}|||\bar L_n\Psi^+\varphi^{(h)}|||^2_{\frac{1}{2}}\lesssim \sum_{j=1}^{n-1}\no{\bar L_j\Psi^+\varphi}_b^2+\no{\Psi^+\varphi}_b^2\quad
\T{for any } \varphi\in C_c^\infty(U\cap b\Om^+).
\end{eqnarray}\end{lemma}

{\it Proof. } We start with \eqref{4.12}. Since $\T{supp }\psi^0\subset \mathcal C^0$, then we have 
\begin{eqnarray}\begin{split}
(1+|\xi'|^2) |\psi^0(\xi)|^2&\lesssim \left(1+\sum_{j=1}^{n-1}|\sigma(L_j)|^2((z_0,0),\xi')\right)|\psi^0(\xi)|^2
\\
&\simleq \left(1+\sum_{j=1}^{n-1}|\sigma(L_j)|^2((x,r),\xi')\right)|\psi^0(\xi)|^2\\
&+\left(\sum_{j=1}^{n-1}|\sigma(L_j)|^2((z_0,0),\xi')-\sum_{j=1}^{n-1}|\sigma(L_j)|^2((x,r),\xi')\right)|\psi^0(\xi)|^2
\\
&\simleq \left(1+\sum_{j=1}^{n-1}|\sigma(L_j)|^2((x,r),\xi)\right)|\psi^0(\xi)|^2\\
&+\T{diam}(\bar\Om^+\cap U)\left(1+\sum_{j=1}^{n-1}|\xi_j|^2\right) |\psi^0(\xi)|^2.
\end{split}
\end{eqnarray}
Hence
\begin{eqnarray}\Label{4.16}
\begin{split}
\no{\Psi^0 \varphi }^2_1\lesssim& \sum_{j=1}^{n}\no{\bar L_j\Psi^0\varphi}^2+\no{\Psi^0\varphi}^2+\T{diam}(U\cap \bar\Om^+) \sum_{j=1}^{2n-1} \no{D_j\Psi^0\varphi }^2.
\end{split}
\end{eqnarray}
 The estimate \eqref{4.12} follows from \eqref{4.16} by taking $U$ sufficiently small so that the last term is absorbed in the left hand side of the estimate.\\

Next, we prove \eqref{4.13}. For all $\varphi\in C_c^\infty(U\cap \bar \Om^+)$,  let $\varphi^{(h)}$ be the harmonic extension of $\varphi_b=\varphi|_{U\cap b\Om^+}$. We have
\begin{eqnarray}\Label{4.17}
|||\Psi^-\varphi|||_1^2
\lesssim \no{\Psi^-(\varphi - \varphi^{(h)})}_1^2+|||\Psi^-\varphi^{(h)}|||^2_1.
\end{eqnarray}

We estimate now $|||\Psi^-\varphi^{(h)}|||_1^2$; we have 
\begin{eqnarray}
\begin{split}
\bar L_n \Psi^-\varphi^{(h)} (x,r)=&(2\pi)^{-2n+1}\int e^{ix\cdot \xi} e^{r (\mu(x,\xi))}\Big(\mu(x,\xi)-\sigma_b(T)(x,\xi)\\
&~~~~~~~~~~~~~~~~~~~~~~~~~~~~~~~~~~~+rp^1(x,\xi)\Big)\psi^-(\xi)\tilde{\varphi}(\xi,0)d\xi  
\end{split}
\end{eqnarray}
where $p^1(x,\xi)$ is the symbol of order 1 and whose associated operator we have denoted by $P^1$.
Choosing $U$ sufficiently small we have  $\sigma(T)_b(x,\xi)\le 0$ when $\xi \in\T{supp}(\psi^-)\subset \mathcal C^-$.  Then, 
$$\mu(x,\xi)-\sigma_b(T)(x,\xi)\simge |\xi|+1.$$
It follows
\begin{eqnarray}\Label{4.19}
\begin{split}
|||\Psi^-\varphi^{(h)}|||_1^2\lesssim& \no{\bar L_n\Psi^-\varphi^{(h)}}^2+\no{rP^1\Psi^-\varphi^{(h)}}^2.\\
\end{split}
\end{eqnarray}
Applying Lemma \ref{l4.10} and inequality \eqref{4.16b} to  the second term in \eqref{4.19}, we get
\begin{eqnarray}\Label{4.20}
\begin{split}
\no{rP_1\Psi^-\varphi^{(h)}}^2 
\lesssim& \no{\Lambda^{-1/2}\Psi^-\varphi}_b^2\\
\lesssim& \no{\Lambda^{-1} D_r\Psi^-\varphi}^2+\no{\Psi^-\varphi}^2\\
\lesssim& \no{\bar L_n\Psi^-\varphi}^2+\no{\Psi^-\varphi}^2.
\end{split}
\end{eqnarray}
For the first term in \eqref{4.19}, we have 
\begin{eqnarray}\Label{4.21}
\begin{split}
\no{\bar L_n\Psi^-\varphi^{(h)}}^2
\lesssim& \no{\bar L_n\Psi^-(\varphi - \varphi^{(h)})}^2+\no{\bar L_n\Psi^-\varphi}^2\\
\lesssim& \no{\Psi^-(\varphi - \varphi^{(h)})}_1^2+\no{\bar L_n\Psi^-\varphi}^2.
\end{split}
\end{eqnarray}\\

Combining \eqref{4.17}, \eqref{4.19}, \eqref{4.20} and \eqref{4.21}, we get
\begin{eqnarray}
\begin{split}
|||\Psi^-\varphi|||_1^2\lesssim \no{\Psi^-(\varphi-\varphi^{(h)})}_1^2+\no{\bar L_n \Psi^-\varphi}^2+\no{\Psi^-\varphi}^2.
\end{split}
\end{eqnarray}
Finally, we estimate $\no{\Psi^-(\varphi-\varphi^{(h)})}_1^2$. Since $\Psi^-(\varphi-\varphi^{(h)})=0$ on $b\Om$, it follows
\begin{eqnarray}\Label{4.23}
\begin{split}
\no{\Psi^-(\varphi-\varphi^{(h)})}^2_1\lesssim&\no{\Delta\Psi^-(\varphi-\varphi^{(h)})}^2_{-1}\\
\lesssim &\no{\Delta\Psi^- \varphi}^2_{-1}+\no{\Delta\Psi^-\varphi^{(h)}}^2_{-1}\\
\lesssim &\sum_{j=1}^{n}\no{ L_j\bar  L_j\Psi^-\varphi}^2_{-1}+\no{P^1\Psi^-\varphi}_{-1}^2+ \no{(rP^2+P^1)\Psi^-\varphi^{(h)}}^2_{-1}\\
\lesssim &\sum_{j=1}^{n}\no{\bar  L_j\Psi^-\varphi}^2+\no{\Psi^-\varphi}^2.\\
\end{split}
\end{eqnarray}
Here the third inequality in \eqref{4.23} follows from \eqref{4.15b}.  This completes the proof of \eqref{4.13}.\\

We prove now \eqref{4.14}. For any $\varphi\in C_c^\infty(U\cap b\Om^+)$, we have 
\begin{equation}
\Label{4.23bis}
\bar L_n \Psi^+\varphi^{(h)} (x,r)=(2\pi)^{-2n+1}\int e^{ix\cdot \xi} e^{r (\mu(x,\xi))}\Big(\mu(x,\xi)-\sigma_b(T)(x,\xi)+rp^1(x,\xi)\Big)\psi^+(\xi)\tilde{\varphi}(\xi,0)d\xi.
\end{equation}

Choosing $U$ sufficiently small we have $\sigma_b(T)(x,\xi)>0$ when $\xi\in \T{supp}\psi^+\subset\mathcal C^+$, so that 
$$-1+\mu-\sigma_b(T)=\sum_{j=1}^{n-1}\frac{\sigma_b( L_j)}{\mu+\sigma_b(T)}\sigma_b(\bar L_j).$$
Since the symbols $\Big\{\dfrac{\sigma_b( L_j)}{\mu+\sigma_b(T)}\Big\}_{1\le j\le n-1}$ are absolutely bounded, $\mu-\sigma_b(T)$ is the symbol of a pseudodifferential operator of the form $\sum_{j=1}^{n-1}P_j\bar L_j+P_0$ where $P_j$ are zero-order operator for $0\le j\le n-1$. We obtain
\begin{eqnarray}
\begin{split}
|||\bar L_n \Psi^+ \varphi^{(h)}|||_{\frac{1}{2}}^2\lesssim &\sum_{j=1}^{n-1}|||(\bar L_j\Psi^+\varphi)^h|||_{\frac{1}{2}}+|||rP_1\Psi^+\varphi^{(h)}|||_{\frac{1}{2}}^2+\sum_{j=1}^{n-1}\no{E_j\Psi^+\varphi}^2_{\frac{1}{2}}\\
\lesssim &\sum_{j=1}^{n-1}||\bar L_j\Psi^+\varphi||_b^2+\no{\Psi^+\varphi}^2_b.
\end{split}
\end{eqnarray}
Here, we used that $\no{E_j\Psi^+\varphi}^2_{\frac{1}{2}}\lesssim \no{\Psi^+\varphi}^2_b$ since $E_j$ is a Poisson operator of order zero.

$\hfill\Box$

Using Lemma~\ref{l4.6} for coefficients of forms, we obtain 

\begin{lemma}\Label{l4.7} 
Let $\Om^+$ be a pseudoconvex at  $z_0$. Then, for a suitable neighborhood $U$ of $z_0$ and for any $u\in C_c^\infty(U\cap \bar\Om^+)^k\cap \T{Dom}(\dib^*)$ with  $k\ge 1$, we have
  $$|||\Psi^0u|||_1^2+|||\Psi^-u|||_1^2\lesssim Q(u,u).$$
Moreover, 
for any $u\in C_c^\infty(U\cap b\Om^+)^k$ with  $k\ge 1$, we have
$$|||\bar L_n \Psi^+ (u^{+})^{(h)}|||_{\frac{1}{2}}^2\lesssim Q_b(u^+,u^+).$$
\end{lemma}

Similarly, we get the basic microlocal estimates for $\Om^-$.
\begin{lemma} Let $\Om^-$ be  pseudoconcave  at  $z_0$. If U is  a sufficiently small  neighborhood of $z_0$, then
we have the three estimates which follow
 \begin{eqnarray}|||\Psi^0\varphi|||^2_{1}\lesssim \sum_{j=1}^{n-1}\no{L_j \Psi^0\varphi}^2+\no{\bar L_n\Psi^0\varphi}^2+\no{\Psi^0\varphi}^2 \quad\T{for any } \varphi \in C_c^\infty(U\cap \bar\Om^-),
 \end{eqnarray}
 \begin{eqnarray}|||\Psi^-\varphi|||^2_{1}\lesssim \sum_{j=1}^{n-1}\no{L_j\Psi^- \varphi}^2+\no{\bar L_n\Psi^0\varphi}^2+\no{\Psi^-\varphi}^2\quad
\T{for any } \varphi \in C_c^\infty(U\cap \bar\Om^-),
\end{eqnarray} \begin{eqnarray}|||\bar L_n\Psi^+\varphi^{(h)}|||^2_{\frac{1}{2}}\lesssim \sum_{j=1}^{n-1}\no{L_j \Psi^+\varphi}_b^2+\no{\Psi^+\varphi}_b^2\quad
\T{for any } \varphi\in C_c^\infty(U\cap b\Om^-).
\end{eqnarray} 
\end{lemma}

\begin{lemma}\Label{l4.8} 
Let $\Om^-$ be pseudoconcave at $z_0$. Then, for a suitable neighborhood $U$ of $z_0$ and for any $u\in C_c^\infty(U\cap \bar\Om^-)^k\cap \T{Dom}(\dib^*)$ with $k\le  n-2$, we have
  $$|||\Psi^0u|||_1^2+|||\Psi^+u|||_1^2\lesssim Q(u,u).$$
Moreover, for any  $u\in C_c^\infty(U\cap b\Om^-)^k$  with $k\le n-2$, we have
$$|||\bar L_n \Psi^- (u^{-})^{(h)}|||_{\frac{1}{2}}^2\lesssim Q_b(u^-,u^-).$$
\end{lemma}

\section{The equivalence of $(f\T-\M)^k$ estimate on $\Om$ and $b\Om$}
In this section, we give the proof  of Theorem \ref{t2}.  This is a consequence of the three theorems which follow, that is, Theorem \ref{t4.10}, \ref{t4.11} and \ref{t4.12}.
\begin{theorem}\Label{t4.10}
Let $\Om^+\subset\C^n$ be a smooth pseudoconvex domain   with boundary $M=b\Om$ at $z_0\in b\Om$. Then,
\begin{enumerate}
  \item[(i)] $(f\T-\M)^k_{\Om^+}$ implies $(f\T-\M_b)^k_{M,+}$  where $\M_b$ is the restriction of  $\M$ to $M$.
  \item[(ii)] $(f\T-\M_b)^k_{M,+}$ implies $(f\T-\M)^k_{\Om+}$ where $\M$ is an extension of $\M$ from $M$ to $\Om$, that is, $\M\big|_M=\M_b$.
\end{enumerate}
for any $k\ge 1$.
\end{theorem}

{\it Proof.} {\it (i). } We need to show that over a neighborhood $U$ of $z_0$ the inequality
$$\no{f(\La)\M_b u^+}_b^2\lesssim Q_b(u^+,u^+)+\no{\Psi^+u}^2_{b,-1}+C_{\M} \no{u^+}^2_{b,-1}$$
holds for any $u\in C^\infty_c(U\cap M)^k$.   Let $\chi=\chi(r)$ be a cut off function  with $\chi(0)=1$.   Applying inequality \eqref{4.16b}, we have
\begin{eqnarray}\Label{4.24}
\begin{split}
\no{f(\Lambda)\M_b u^+}^2_b\lesssim& \no{\Lambda^{\frac{1}{2}} f(\Lambda)\M\chi u^{+\,(h)}}^2 +\no{\Lambda^{-\frac{1}{2}}  D_r(  f(\Lambda)\M\chi u^{+\,(h)})}^2\\
\lesssim&\no{f(\Lambda)\M\chi\zeta'  R^+  u^{+\,(h)}}^2\\
& +\no{\Lambda^{-1} f(\Lambda) D_r(\M\chi\zeta'R^+  u^{+\,(h)})}^2+error,
\end{split}
\end{eqnarray}
where $\zeta'=1$ on supp$(u^+)$ and $\T{supp}(\chi\zeta')\subset\subset U'$. Here the second inequality follows the same argument in the proof of Lemma~\ref{l4.5} and therefore  the error term is estimated by
\begin{eqnarray}
\begin{split}error\lesssim& \no{\La^{\frac{1}{2}} u^{+\,(h)}}^2+\no{\La^{-\frac{1}{2}}D_r u^{+\,(h)}}^2+C_\M\left(\no{\La^{-\frac{1}{2}} u^{+\,(h)}}^2+\no{\La^{-\frac{3}{2}}D_r u^{+\,(h)}}^2\right)\\
&+\no{\La^{-\frac{1}{2}}\chi\zeta'\Psi^+u^{(h)}}^2+\no{\La^{-\frac{3}{2}}D_r \chi\zeta'\Psi^+u^{(h)}}^2\\
\lesssim&\no{ u^{+}}_b^2+C_\N\no{u^+}_{b,-1}^2+\no{\Psi^+u}_{b,-1}^2
\end{split}
\end{eqnarray}
where the last inequality follows Lemma \ref{l4.10}.
Notice that $\chi\zeta'  R^+  u^{+\,(h)}\in \T{Dom}(\dib^*)$. Using the hypothesis of the theorem to estimate the second  line of \eqref{4.24} and applying Lemma \ref{l3.38} to the second term in the last line of \eqref{4.24}, we have that \eqref{4.24} can be continued by
\begin{eqnarray}
\begin{split}
\lesssim& Q(\chi\zeta'  R^+  u^{+\,(h)}, \chi\zeta'  R^+  u^{+\,(h)})+C_\M\no{\chi\zeta'  R^+  u^{+\,(h)}}_{-1}^2+\no{ u^{+}}_b^2+C_\M\no{u^+}_{b,-1}^2+\no{\Psi^+u}_{b,-1}^2\\
\lesssim& \sumK\sum_{ij}^{n-1}(r_{ij}\zeta'  R^+  u^{+}_{iK}, \zeta'  R^+  u^{+}_{jK})_b\\
&+\sumJ\sum_{j=1}^{n}\no{\bar L_j\chi\zeta'  R^+  u_J^{(h)+}}^2+\no{\chi\zeta'  R^+  u^{+\,(h)}}^2+C_\M\no{\chi\zeta'  R^+  u^{+\,(h)}}_{-1}^2\\
&+\no{ u^{+}}_b^2+C_\M\no{u^+}_{b,-1}^2+\no{\Psi^+u}_{b,-1}^2\\
\lesssim& Q_b(u^+,u^+)+C_\M\no{u^{+}}_{b,-1}^2+\no{\Psi^+u}^2_{b,-1}\\
&+\sumJ\sum_{j=1}^{n-1}\no{\Lambda^{\frac{1}{2}} (\bar L_j u_J^{+})^{(h)}}^2+\no{\La^{\frac{1}{2}} \bar L_n \Psi^+(u^{+})^{(h)}}^2\\
\lesssim& Q_b(u^+,u^+)+C_\M\no{u^{+}}_{b,-1}^2+\no{\Psi^+u}^2_{b,-1},
\end{split}
\end{eqnarray}
where the second inequality follows from \eqref{new1} (with the choice $q_0=0$), the third from Lemma \ref{l4.5} and the last from Lemma \ref{l4.7}.\\

{\it  (ii). }  For any $u\in C^\infty_c(U\cap \bar\Om)^k\cap \T{Dom}(\dib^*)$, we decompose $u=u^\tau+u^\nu$ and  $u^\tau=u^{\tau\,+   }+u^{\tau\,-   }+u^{\tau\,0   }$. Since $u^\nu$ satisfies elliptic estimates and on account of  Lemma \ref{l4.7}, we have
\begin{eqnarray}
\begin{split}
\no{f(\Lambda)\M u^\nu}^2\le |||u^\nu|||_1^2+C_\M\no{u^\nu}_{-1}^2&\lesssim Q(u,u)+C_\M\no{u}_{-1}^2,\\
\no{f(\Lambda)\M u^{\tau\,0   }}^2\le |||u^{\tau\,0   }|||_1^2+C_\M\no{u^{\tau\,0   }}_{-1}^2&\lesssim Q(u,u)+C_\M\no{u}_{-1}^2,\\  
\no{f(\Lambda)\M u^{\tau\,-   }}^2\le |||u^{\tau\,-   }|||_1^2+C_\M\no{u^{\tau\,-   }}_{-1}^2&\lesssim Q(u,u)+C_\M\no{u}_{-1}^2.\\    
\end{split}
\end{eqnarray}
Moreover, using Lemma 1.8. and Proposition 1.9 in \cite{KZ10}, we have
\begin{eqnarray}
\begin{split}
\no{f(\Lambda)\M u^{\tau\,+   }}^2\lesssim& \no{\Lambda^{-\frac{1}{2}}f(\Lambda)\M_b u^{\tau\,+   }_b}^2_b+Q(u^\tau,u^\tau)+C_\M\no{u^\tau}_{-1}^2\\    
\lesssim& \no{\Lambda^{-\frac{1}{2}}f(\Lambda)\M_b u^{\tau\,+   }_b}^2_b+Q(u,u)+C_\M\no{u}_{-1}^2.\\    
\end{split}
\end{eqnarray}
Thus, we obtain
\begin{eqnarray}
\begin{split}
\no{f(\Lambda)\M u}^2\lesssim &\no{f(\Lambda)\M u^{\tau\,+   }}^2+\no{f(\Lambda)\M u^{\tau\,-   }}^2+\no{f(\Lambda)\M u^{\tau\,0   }}^2+\no{f(\Lambda)\M u^\nu}^2\\
\lesssim& \no{\Lambda^{-\frac{1}{2}}f(\Lambda)\M_b u^{\tau\,+   }_b}^2_b+Q(u,u)+C_\M\no{u}_{-1}^2.\\    
\end{split}
\end{eqnarray}
Hence,  we only need to estimate $\no{\Lambda^{-\frac{1}{2}}f(\Lambda)\M_b u^{\tau\,+   }_b}^2_b$.
We begin by noticing that
\begin{eqnarray}\Label{4.30}
\begin{split}
\no{\Lambda^{-\frac{1}{2}}f(\Lambda)\M_b u^{\tau\,+   }_b}^2_b
\lesssim& \no{f(\Lambda)\M_b (\zeta'\Lambda^{-\frac{1}{2}}u^{\tau})^+_b}^2_b+C_\M\no{u^{\tau\,}_b}^2_{b,-1}\\
\end{split}
\end{eqnarray}
where $\zeta'\equiv1$ over $\T{supp}u^{\tau\,+}$.
Using the hypothesis we can continue  \eqref{4.30} by
\begin{eqnarray}
\begin{split}
\lesssim&Q_b((\zeta'\Lambda^{-\frac{1}{2}}u^{\tau})^+_b,(\zeta'\Lambda^{-\frac{1}{2}}u^{\tau})^+_b)+\no{\Psi^+(\zeta'\Lambda^{-\frac{1}{2}}u^{\tau})_b}_{b,-1}^2\\
&+C_\M\no{(\zeta'\Lambda^{-\frac{1}{2}}u^{\tau})^+_b}^2_{b,-1}+C_\M\no{u^{\tau }_b}^2_{b,-1}\\
\lesssim&\sumK\sum_{ij=1}^{n-1}(r_{ij}\zeta'''R^+(\zeta'\Lambda^{-\frac{1}{2}}u^{\tau})^+_{b,iK},\zeta'''R^+(\zeta'\Lambda^{-\frac{1}{2}}u^{\tau})^+_{b,jK})_b\\
&\sumJ\sum_{j=1}^{n-1}\no{\bar L_j( \zeta'\Lambda^{-\frac{1}{2}}u^{\tau})^+_{b,J}}_b^2+\no{(\zeta'\Lambda^{-\frac{1}{2}}u^{\tau})^+_b}^2_b\\
&+\no{\Psi^+(\zeta'\Lambda^{-\frac{1}{2}}u^{\tau})_b}_{b,-1}^2+C_\M\no{u^{\tau }_b}^2_{b,-1}\\
\lesssim&Q(\zeta'''R^+\zeta''\Psi^+\zeta'\Lambda^{-\frac{1}{2}}u^{\tau},\zeta''R^+\zeta''\Psi^+\zeta'\Lambda^{-\frac{1}{2}}u^{\tau})\\
&+\sumJ\sum_{j=1}^{n-1}\no{\bar L_j( \zeta'\Lambda^{-\frac{1}{2}}u^{\tau})^+_{b,J}}_b^2+\no{(\zeta'\Lambda^{-\frac{1}{2}}u^{\tau})^+_b}^2_b+C_\M\no{u^{\tau }_b}^2_{b,-1}\\
\end{split}
\end{eqnarray}
Since $\zeta'''R^+\zeta''\Psi^+\zeta'\Lambda^{-\frac{1}{2}}$ is a tangential pseudodifferential operators of order zero, then 
 $$Q(\zeta'''R^+\zeta''\Psi^+\zeta'\Lambda^{-\frac{1}{2}}u^{\tau},\zeta'''R^+\zeta''\Psi^+\zeta'\Lambda^{-\frac{1}{2}}u^{\tau})\lesssim Q(u^{\tau},u^{\tau}).$$

To estimate the last line  we proceed as follows. Since $\bar L_j( \zeta'\Lambda^{-\frac{1}{2}}u^{\tau})^+\in C^\infty_c(U\cap \bar\Om)$, then using inequality \eqref{4.16b}, we have
\begin{eqnarray}
\begin{split}
\no{\bar L_j( \zeta'\Lambda^{-\frac{1}{2}}u^{\tau})^+_{b}}_b^2\lesssim&\no{\Lambda^{\frac{1}{2}}\bar L_j( \zeta'\Lambda^{-\frac{1}{2}}u^{\tau})^+}^2+\no{\Lambda^{-\frac{1}{2}}D_r\bar L_j( \zeta'\Lambda^{-\frac{1}{2}}u^{\tau})^+}^2\\
\lesssim&\no{\bar L_j u^{\tau}}^2+\no{\Lambda^{-1}D_r\bar L_j u^{\tau}}^2+\no{u^{\tau}}^2\\
\lesssim&\no{\bar L_j u^{\tau}}^2+\no{T\Lambda^{-1}\bar L_j u^{\tau}}^2+\no{\bar L_n\Lambda^{-1}\bar L_j u^{\tau}}^2+\no{u^{\tau}}^2\\
\lesssim&\no{\bar L_j u^{\tau}}^2+\no{\bar L_n u^{\tau}}^2+\no{u^{\tau}}^2\\
\lesssim&Q(u^\tau,u^\tau),
\end{split}
\end{eqnarray}
and similarly $\no{(\zeta'\Lambda^{-\frac{1}{2}}u^{\tau})^+_b}^2_b\lesssim Q(u^\tau,u^\tau)$. We finish this proof with the estimate of  $C_\M\no{u^{\tau}_b}_{b,-1}^2$. Using the interpolation inequality  
$$C_\M\no{u^{\tau}_b}_{b,-1}^2\lesssim \no{D_r u^{\tau}}^2_{-1}+C_\M\no{u^{\tau}}^2_{-1}\lesssim Q(u^\tau,u^\tau)+C_\M\no{u^{\tau}}^2_{-1}$$
This concludes the proof of Theorem \ref{t4.10}.

$\hfill\Box$

Similarly, we get the equivalence of $(f\T-\M)^k$ on $\Om^-$ and $M$

\begin{theorem}\Label{t4.11}
Let $\Om^-$ be a smooth pseudoconcave domain at $z_0\in b\Om$. 
Then $(f\T-\M)^k_{\Om^-}$ is equivalent to $(f\T-\M_b)^k_{b\Om,-}$ for $\M|_{b\Om}=\M_b$ for any $k\le n-2$.
\end{theorem}

\section{The equivalence of microlocal estimates $b\Om$}
In this section, we prove the equivalence of microlocal estimates on hypersurface.
\begin{theorem}\Label{t4.12}
Let $M$ be a hypersurface (not necessarily pseudoconvex) and $z_0\in M$.  Then $(f\T-\M)^k_{M,+}$ holds at $z_0$ if and only if $(f\T-\M)^{n-1-k}_{M,-}$ holds at $z_0$.
\end{theorem}

{\it Proof.} We define the local conjugate-linear duality map $F^k : \A^{0,k}_b\to \A^{0,n-1-k}_b$ as follows. If $u=\sumJ u_J\bom_J$ then 
$$F^{k}u=\underset{|J'|=n-1-k}{{\sum}'} \epsilon^{\{J ,J'\}}_{\{1,..., n-1\}} \bar u_J\bom_{J'} , $$
where $J'$ denotes the strictly increasing $(n-k-1)$-tuple consisting of all integers in $[1,n-1]$ which do not belong to $J$ and $\epsilon^{J,J'}_{\{1,n-1\}}$ is the sign of the permutation $\{J, J'\}\simto \{1,\dots, n-1\}$. \\

By this definition, we obtain $F^{n-1-k}F^ku=u$, $\no{F^ku}_b=\no{u}_b$, $\dib_b F^k u=F^{k-1}\dib^*_b u +\cdots$, and 
$\dib^*_b F^k u=F^{k+1}\dib_b u +\cdots,$
for any $u\in C^\infty_c(U\cap M)^k$, where dots refers the term in which $u$ is not differentiated. We get
\begin{eqnarray}\Label{a1} Q_b(F^ku,F^ku)\cong Q_b(u,u).
\end{eqnarray}

We consider two cases of multiplier.\\
{\bf Case 1.} If $\M$ is a function, then 
\begin{eqnarray}\Label{a2}
\overline{\M} F^ku=F^k\left(\M u\right).
\end{eqnarray}
{\bf Case 2.} Let $\M=\underset{j=1}{\overset{n-1}{\sum}}\M_j\om_j\in \A_b^{1,0}$. 
We define the operator $\overline \M: \A_b^{k}\to \A_b^{k+1}$ by $\overline{\M} u:=\sumJ\underset{j=1}{\overset{n-1}{\sum}}\overline\M_j u_J\bom_j\we\bom_J;$ and  $\overline{\M}^*:\A_b^{k}\to\A_b^{k-1}$ by $\overline{\M}^* u:=\sumK\underset{j=1}{\overset{n-1}{\sum}}\M_j u_{jK}\bom_K;$ then we obtain 
\begin{eqnarray}\Label{a3}\overline{\M}F^k u=F^{k-1}(\overline{\M}^* u).
\end{eqnarray}
 
We notice that with the definitions of $\overline{\M}$ and $\overline{\M}^*$, above coincide with ones in \eqref{Mpm}, i.e., that $\M u^+=\overline{\M}^*u^+$ and $\M u^-=\overline{\M}u^-$. Replace $u$ by $u^+$  and $u^-$ in \eqref{a1}, \eqref{a2} and \eqref{a3}, we obtain 
\begin{eqnarray}\begin{cases} Q_b(F^ku^+,F^ku^+)&\cong Q_b(u^+,u^+)\\
\no{f(\La)\M F^ku^+}_b&=\no{f(\La)\M u^+}_b.
\end{cases}
\end{eqnarray}
and 
\begin{eqnarray}\begin{cases} Q_b(F^ku^-,F^ku^-)&\cong Q_b(u^-,u^-)\\
\no{f(\La)\M F^ku^-}_b&=\no{f(\La)\M u^-}_b.
\end{cases}
\end{eqnarray}

On the other hand, we have 
\begin{eqnarray}
\begin{split}
\overline{(u^+)}(x)=&(2\pi)^{2n-1}\overline{\int_{\R^{2n-1}_\xi}e^{ix\xi} \psi^+(\xi) \int_{\R^{2n-1}_y}e^{-iy\xi}u(y)dyd\xi}\\
=&(2\pi)^{2n-1}\int_{\R^{2n-1}_\xi}e^{-ix\xi} \psi^+(\xi) \int_{\R^{2n-1}_y}e^{iy\xi}\bar{u}(y)dyd\xi\\
\overset{\xi:=-\xi}=&(-2\pi)^{2n-1}\int_{\R^{2n-1}_\xi}e^{ix\xi} \psi^+(-\xi) \int_{\R^{2n-1}_y}e^{-iy\xi}\bar{u}(y)dyd\xi\\
=&(-2\pi)^{2n-1}\int_{\R^{2n-1}_\xi}e^{ix\xi} \psi^-(\xi) \int_{\R^{2n-1}_y}e^{-iy\xi}\bar{u}(y)dyd\xi\\
=&-(\overline{u})^-(x),
\end{split}
\end{eqnarray}
for any $u\in C^\infty_c(U\cap M)$. Hence, 
$$F^{k}u^+=-\sum \epsilon^{\{J ,J'\}}_{\{1,..., n-1\}} (\bar u)^-_J\bom_{J'}, $$
and 
$$F^{k}u^-=-\sum \epsilon^{\{J ,J'\}}_{\{1,..., n-1\}} (\bar u)^+_J\bom_{J'}, $$
for any $u\in C^\infty_c(U\cap M)^k$. This proves Theorem \ref{t4.12}.

$\hfill\Box$

It is also interesting to remark that when $M$ is pseudoconvex hypersurface, then $(f\T-\M)_{M,+}^k$ implies $(f\T-\M)_{M,+}^{k+1}$, and $(f\T-\M)_{M,-}^k$ implies $(f\T-\M)_{M,-}^{k-1}$.
\begin{lemma}\Label{l4.12}
Let $M$ be a pseudoconvex hypersurface. 
\begin{enumerate} 
  \item Assume that $(f\T-\M)_{M,+}^{k_o}$ holds then  $(f\T-\M)_{M,+}^{k}$ also holds for any $k_o\le k\le n-1$.
\item Assume that $(f\T-\M)_{M,-}^{k_o}$ holds then  $(f\T-\M)_{M,-}^{k}$ also holds for any $0\le k\le k_o$.
\end{enumerate}
\end{lemma}
The proof is analogous to Lemma 3.11 in \cite{Kha10} or using the combination of that lemma with Theorem  \ref{t4.10} and Theorem \ref{t4.11}.

\begin{corollary}\Label{c4.12}
Let $M$ be a pseudoconvex hypersurface at $z_0$, let $\M_b \in \A^{0,0}_b$ and $\M \in \A^{0,0}$ such that $\M|_M=\M_b$.
 Then,  for $1\le k\le n-2$, the following are equivalent:
\begin{enumerate}
  \item $(f\T-\M_b)^k_{M}$ holds.
  \item Both $(f\T-\M_b)^k_{M,+}$ and $(f\T-\M_b)^k_{M,-}$ hold.
%\item Both $(f\T-\N)^k_{M,+}$ and $(f\T-\N)^{n-1-k}_{M,+}$ hold. 
  \item $(f\T-\M_b)^l_{M,+}$  holds with $l=\min\{k,n-1-k\}$.
  \item $(f\T-\M_b)^l_{M,-}$  holds with $l= \max\{k,n-1-k\}$
  \item Both $(f\T-\M)^k_{\Om^+}$ and $(f\T-\M)^k_{\Om^-}$ hold.
  %\item Both $(f\T-\M)^k_{\Om^+}$ and $(f\T-\M)^{n-1-k}_{\Om^+}$ hold. 
  \item $(f\T-\M)^l_{\Om^+}$ holds with $l= \min(k,n-1-k)$.
  \item $(f\T-\M)^l_{\Om^-}$ holds with  $l= \max(k,n-1-k)$.
\end{enumerate}
\end{corollary}

{\it Proof. }
The proof follows from Theorem~\ref{t4.10}, \ref{t4.11}, \ref{t4.12} and Lemma \ref{l4.12} combined with the fact that $\no{u^0}^2_{b,1}\simleq Q_b(u,u)$.

$\hfill\Box$

\section*{Acknowledgments} 
This research is partially supported by  Vietnam National Foundation for Science and Technology Development (NAFOSTED) under grant number 101.01-2012.16. The author gratefully acknowledge the careful reading by the referee. The exposition and rigor of the paper were improved by the close reading.

\bibliographystyle{alpha}
%\bibliography{Khanh}

\end{document}